\input amstex
\magnification=\magstep1
\def\ds{\baselineskip 20pt plus 2pt}
\def\ss{\baselineskip 15pt plus 1.5pt}
\ss
\def\Bbb#1{\hbox{\bf #1}}
\def\N{\Bbb N}
\def\Z{\Bbb Z}
\def\C{\Bbb C}
\def\R{\Bbb R}
\def\supp{\operatorname{supp}}
\def\It{\operatorname{int}}
\voffset=0.4in
\font\bigbold=cmbx10 scaled \magstep2
\font\sl=cmbsy10
\def\slK{\hbox{$K\textfont1=\sl$}}
\def\sqr{\vcenter {\hrule height.3mm
\hbox {\vrule width.3mm height 2mm \kern2mm
\vrule width.3mm } \hrule height.3mm }}
\def\endproof{\nobreak\hfill$\sqr$\bigskip\goodbreak}

\pageno=0
\footline={\ifnum\pageno=0\hfill\else\hss\tenrm\folio\hss\fi}
\def\ans{\vrule height.1pt width80pt depth0pt}
\centerline {\bigbold p-Summing Operators on}

\centerline {\bigbold Injective Tensor Products of Spaces}
\vskip 1truein
\centerline {by}
\vskip 1truein
\centerline {\bf Stephen Montgomery-Smith$^{(*)}$ and Paulette 
Saab$^{(**)}$}
\vskip1truein
{\narrower\smallskip\noindent {\bf Abstract}\ \ Let $X,Y$ and $Z$ be
Banach spaces, and let $\prod_p(Y,Z)\ (1\leq p<\infty)$ denote the
space of $p$-summing operators from $Y$ to $Z$.  We show that, if 
$X$ is a {\it \$}$_\infty$-space, then a bounded linear operator $T:\
X\hat \otimes_\epsilon Y\longrightarrow Z$ is 1-summing if and only if 
a
naturally associated operator $T^\#:\ X\longrightarrow \prod_1(Y,Z)$ 
is
1-summing.  This result need not be true if $X$ is not a {\it
\$}$_\infty$-space.  For $p>1$, several examples are given with
$X=C[0,1]$ to show that $T^\#$ can be $p$-summing without $T$ being
$p$-summing.  Indeed, there is an operator $T$ on
$C[0,1]\hat \otimes_\epsilon \ell_1$ whose associated operator $T^\#$ 
is
2-summing, but for all $N\in \N$, there exists an $N$-dimensional
subspace $U$ of $C[0,1]\hat \otimes_\epsilon \ell_1$ such that $T$
restricted to $U$ is equivalent to the identity operator on
$\ell^N_\infty$. Finally, we show that there is a compact Hausdorff 
space $K$\
and a bounded linear operator $T:\ C(K)\hat \otimes_\epsilon
\ell_1\longrightarrow \ell_2$ for which $T^\#:\ C(K)\longrightarrow
\prod_1(\ell_1, \ell_2)$ is not 2-summing.

\smallskip}
\vfill

\ans

\item{$^{(*)}$} Research supported in part by an NSF Grant DMS 9001796
\item{$^{(**)}$} Research supported in part by an NSF Grant DMS 
87500750
\par\noindent A.M.S.\ (1980) subject classification: 46B99

\eject
\voffset=-.3in
\ds

\noindent {\bf Introduction}\ \ Let $X$ and $Y$ be Banach spaces, and
let $X\hat \otimes_\epsilon Y$ denote their injective tensor
product. In this paper, we shall study the behavior of those operators
on
$X\hat \otimes_\epsilon Y$ that are $p$-summing. 

\medskip

If $X$, $Y$ and $Z$ are Banach
spaces, then every $p$-summing operator $T:\ X\hat \otimes_\epsilon
Y\longrightarrow Z$ induces a $p$-summing linear operator $T^\#:\
X\longrightarrow \prod_p(Y,Z)$. This raises the
following question:\ \ given two Banach spaces $Y$ and $Z$, and $1\leq 
p<\infty$,
for what Banach spaces
$X$ is it true that a bounded linear operator $T:\ X\hat
\otimes_\epsilon Y\longrightarrow Z$ is $p$-summing whenever $T^\#:\
X\longrightarrow \prod_p(Y,Z)$ is $p$-summing?  

\medskip

In [11], it was shown that
whenever $X=C(\Omega)$ is a space of all continuous functions on a 
compact
Hausdorff space $\Omega$, then $T:\ C(\Omega)\hat \otimes_\epsilon
Y\longrightarrow Z$ is 1-summing if and only if $T^\#:\ 
C(\Omega)\longrightarrow
\prod_1(Y,Z)$ is 1-summing.
We will extend this result by showing that this result still remains 
true if $X$
is any {\it \$}$_\infty$-space. We will also give an example to show 
that the
result need not be true if $X$ is not a {\it \$}$_\infty$-space.  For 
this,
we shall exhibit a 2-summing operator $T$ on $\ell_2\hat 
\otimes_\epsilon \ell_2$
that is not 1-summing, but such that the associated operator $T^\#$ is 
1-summing.

\medskip

The case $p > 1$ turns out to be
quite different.  Here, the {\it \$}$_\infty$-spaces do not seem to 
play any
important role.  We show that for each $1<p<\infty$, there exists a 
bounded
linear operator $T:\ C[0,1]\hat \otimes_\epsilon \ell_2\longrightarrow 
\ell_2$
such that $T^\#:C[0,1]\longrightarrow \prod_p(\ell_2, \ell_2)$ is 
$p$-summing,
but such that $T$ is not $p$-summing. We
will also give an example that shows that, in general, the condition 
on $T^\#$ to
be 2-summing is too weak to imply any good properties for the operator 
$T$\ at
all.  To illustrate this, we shall exhibit a bounded linear operator 
$T$ on
$C[0,1]\hat \otimes_\epsilon \ell_1$ with values in a certain Banach 
space $Z$,
such that $T^\#:\ C[0,1]\longrightarrow \prod_2(\ell_1,Z)$ is 
2-summing, but for
any given $N\in \N$, there exists a subspace $U$ of $C[0,1]\hat 
\otimes_\epsilon
\ell_1$, with $\dim U=N$, such that $T$ restricted to $U$ is 
equivalent to
the identity operator on $\ell^N_\infty$.

\medskip

Finally, we show that there is a compact Hausdorff space $K$\
and a bounded linear operator $T:\ C(K)\hat \otimes_\epsilon
\ell_1\longrightarrow \ell_2$ for which $T^\#:\ C(K)\longrightarrow
\prod_1(\ell_1, \ell_2)$ is not 2-summing.

\bigskip
\noindent {\bf I - Definitions and Preliminaries}

Let $E$ and $F$ be Banach spaces, and let $1\leq q\leq p<\infty$.  An
operator $T:\ E\longrightarrow F$ is said to be $(p,q)$-{\bf summing} 
if
there exists a constant $C\ge 0$ such that for any finite sequence
$e_1, e_2,\ldots, e_n$ in $E$, we have
$$
\left(\sum\limits^n_{i=1}\parallel T(e_i)\parallel^p\right)^{\frac 
1p}\leq C\sup
\left\{\left(\sum\limits^n_{i=1}|e^*(e_i)|^q\right)^{\frac 1q}:e^*\in
E^*, \parallel e^*\parallel\leq 1\right\} .
$$
We let $\pi_{p,q}(T)$ denote the smallest constant $C$ such that the
above inequality holds, and let $\prod_{p,q}(E,F)$\ be the space of 
all
$(p,q)$-summing operators from $E$ to $F$ with the norm
$\pi_{p,q}$.  It is easy to check that $\prod_{p,q}(E,F)$ is a Banach
space.  In the case $p=q$, we will simply write $\prod_p(E,F)$ and 
$\pi_p$.
We will use the fact that $T\in \prod_{p,q}(E,F)$ if and only if 
$\sum\limits_n
\left\|
Te_n\right\|^{p} <\infty$\ for every
infinite sequence $(e_n)$ in $E$ with $\sum\limits_n
|e^*(e_n)|^q<\infty\text { for each }e^*\in E^*$. That is to say, $T$\ 
is in
$\prod_{p,q} (E,F)$ if and only if
$T$
sends all weakly $\ell_q$-summable sequences into strongly 
$\ell_p$-summable
sequences. In what follows we shall
mainly be interested in the case where $p=q$ and $p=1\text { or }2$.

\medskip

Given two Banach spaces $E$ and $F$, we will let $E\hat
\otimes_\epsilon F$ denote their injective tensor product,
that is, the completion of the algebraic tensor product $E\otimes F$
under the cross norm $\parallel \cdot \parallel_\epsilon$\ given by 
the
following formula. If $\sum\limits^n_{i=1} e_i\otimes x_i\in E\otimes 
F$,
then $$
\parallel \sum\limits^n_{i=1}e_i\otimes x_i\parallel_\epsilon
=\sup\left\{\Bigg|\sum\limits^n_{i=1}e^*(e_i)x^*(x_i)\Bigg|:\ 
\parallel
e^*\parallel \leq 1, \parallel x^*\parallel \leq 1, e^*\in E^*, x^*\in
F^*\right\}.
$$

\medskip

We will say that a bounded linear operator $T$ between two Banach 
spaces $E$
and $F$ is called an {\bf integral operator} if the bilinear form 
$\tau$\ defines
an element of $(E\hat \otimes_\epsilon F^*)^*$, where $\tau$\ is 
induced by
$T$\ according to the formula $\tau (e,x^*)=x^*(Te)$ ($e\in E$,
$x^*\in F^*$). We will define the {\bf integral norm} of $T$,
denoted by $\parallel T\parallel_{\text {int}}$, by
$$
\parallel T\parallel_{\text {int }}
=\sup\left\{\Bigg|\sum\limits^n_{i=1}x_i^*(Te_i)\Bigg|:\ \parallel
\sum\limits^n_{i=1}e_i\otimes x_i^*\parallel_\epsilon \leq 1\right\}.
$$
The space of all integral operators from a Banach space $E$ into a
Banach space $F$ will be denoted by $I(E,F)$.  We note that
$I(E,F)$ is a Banach space under the integral norm
$\parallel\quad\parallel_{\text {int}}$.  

\medskip

We will say that a Banach
space $X$ is a {\it \$}$_\infty$-{\bf space} if, for some $\lambda>1$, 
we have
that for every finite dimensional subspace $B$ of $X$, there exists a 
finite
dimensional subspace $E$ of $X$ containing $B$, and an invertible 
bounded
linear operator $T:\ E\longrightarrow \ell^{\dim E}_\infty$ such that
$\parallel T\parallel\ \parallel T^{-1}\parallel\leq \lambda$.

\medskip

It is well known that for any
Banach spaces $E$ and $F$, if $T$ is in $I(E,F)$, then it is also in
$\prod_1(E,F)$, with $\pi_1(T)\leq \parallel T\parallel_{\text 
{int}}$.  But
$I(E,F)$ is strictly included in $\prod_1(E,F)$.    It was shown
in [12, p. 477] that a Banach space $E$ is a {\it \$}$_\infty$-space 
if and only
if for any Banach space $F$, we have that $I(E,F)=\prod_1(E,F)$. We 
will use
this characterization of {\it \$}$_\infty$-spaces in the sequel.

\medskip

Finally, we note the following characterization of 1-summing operators
(called right semi-integral by Grothendieck in [5]), which will be 
used
later.

\medskip

\noindent {\bf Proposition 1}\ \ Let $E$ and $F$ be Banach spaces. 
Then
the following properties about a bounded linear operator $T$ from $E$ 
to
$F$ are equivalent:

\item {(i)} $T$ is 1-summing;

\item {(ii)} There exists a Banach space $F_1$, and an isometric
injection $\varphi:\ F\longrightarrow F_1$, such that $\varphi\circ 
T:\
E\longrightarrow F_1$ is an integral operator.

\medskip

For all other undefined notions we shall refer the reader to either 
[3],
[7] or [10].

\vfill\eject
\noindent {\bf II\ 1-Summing and Integral Operators}

Let $X$ and $Y$ be Banach spaces with injective tensor
product $X\hat \otimes_\epsilon Y$.  For
a Banach space $Z$, any bounded linear operator $T:\ X\hat
\otimes_\epsilon Y\longrightarrow Z$ induces a linear
operator $T^\#$ on $X$ by
$$
T^\#x(y)=T(x\otimes y) \qquad (y\in Y) .
$$
It is clear that the range of $T^\#$ is the space {\it \$}$(Y,Z)$ of
bounded linear operators from $Y$ into $Z$, and that $T^\#$ is a 
bounded
linear operator.

\medskip

In this section, we are going to investigate the 1-summing operators, 
and the
integral operators, on $X\hat \otimes_\epsilon Y$. We will use 
Proposition~1 to
relate these two ideas together. First of all, we have the following 
result.

\medskip

\noindent {\bf Theorem 2}\ \ Let $X,Y$ and $Z$ be Banach spaces, and 
let
$T:\ X\hat \otimes_\epsilon Y\longrightarrow Z$ be a bounded linear
operator.  Denote by $i:\ Z\longrightarrow Z^{**}$ the isometric
embedding of $Z$ into $Z^{**}$. Then the following two properties are
equivalent:

\item {(i)} $T\in I (X\hat \otimes_\epsilon Y,Z)$;

\item {(ii)} $\hat i \circ T\in I (X, I(Y,Z^{**}))$, where $\hat i:\
I(Y,Z)\longrightarrow I(Y,Z^{**})$ is defined by $\hat i (U)=i\circ U$ 
for
each $U\in I(Y,Z)$.

\noindent In particular, if $T^\#\in I(X, I(Y,Z))$, then $T\in I(X\hat
\otimes_\epsilon Y,Z)$.
\medskip
{\bf Proof:}\ \ First, we show that $(X\hat \otimes_\epsilon Y)\hat
\otimes_\epsilon Z^*$\ and $X\hat \otimes_\epsilon
(Y\hat \otimes_\epsilon Z^*)$\ are isometrically isomorphic to one 
another. To
see this, note that the algebraic tensor product is an associative 
operation,
that is, $(X\otimes Y)\otimes Z^*$ and $X\otimes (Y\otimes Z^*)$ are
algebraically isomorphic. Also, they are both generated by elements of 
the form
$\sum\limits^n_{i=1}x_i\otimes y_i\otimes z_i^*$, where $x_i\in X,\ 
y_i\in Y$
and $z^*_i\in Z^*$. Now, if we let $B(X^*),\ B(Y^*)$ and $B(Z^{**})$ 
denote the
dual unit balls of $X^*,\ Y^*$ and $Z^{**}$ equipped with their 
respective
weak$^{*}$ topologies, then the spaces $(X\otimes_\epsilon 
Y)\otimes_\epsilon
Z^*$\ and $X\otimes_\epsilon (Y\otimes_\epsilon Z^*)$\ embed 
isometrically into
$C\left(B(X^*)\times B(Y^*)\times B(Z^{**})\right)$ in a natural way, 
by
$$
\langle \sum\limits^n_{i=1}x_i\otimes y_i\otimes z_i^*,\quad (x^*, 
y^*,
z^{**})\rangle=
\sum\limits^n_{i=1}x^*(x_i)y^*(y_i)z^{**}(z^*_i) ,
$$
where $\sum\limits^n_{i=1}x_i\otimes y_i\otimes z^*_i$ is in
$(X\otimes_\epsilon Y)\otimes_\epsilon Z^*$\ or 
$X\otimes_\epsilon (Y\otimes_\epsilon Z^*)$, and $(x^*, y^*,
z^{**})$ is in the compact set $B(X^*)\times B(Y^*)\times B(Z^{**})$. 
Thus
both spaces
$(X\hat \otimes_\epsilon Y)\hat \otimes_\epsilon Z^*$ and $X\hat
\otimes_\epsilon (Y\hat \otimes_\epsilon Z^*)$ can be thought of as 
the
closure in $C\left(B(X^*)\times B(Y^*)\times B(Z^{**})\right)$ of the
algebraic tensor product of $X$, $Y$ and $Z^*$.

\medskip

Now let us assume that $T:\ X\hat \otimes_\epsilon Y\longrightarrow
Z$ is an integral operator.  Then the bilinear map $\tau$
on $X\hat\otimes_\epsilon Y\times Z^*$, given by $\tau
(u,z^*)=z^*(Tu)$ for $u\in X\hat \otimes_\epsilon Y$ and $z^*\in Z^*$,
defines an element of $\left(X\hat \otimes_\epsilon Y\hat
\otimes_\epsilon Z^*\right)^*$, that is,
$$
\parallel 
T\parallel_{\It}=\sup\left\{|\sum\limits^n_{i=1}z_i^*\left(T(x_i\otimes
y_i)\right):\ \ \parallel\sum\limits^n_{i=1}x_i\otimes y_i \otimes
z^*_i\parallel_\epsilon \leq 1\right\} .\tag *
$$
To show that for every $x$\ in $X$\ the operator $T^\#x$\ is in 
$I(Y,Z)$, with
$$
\parallel T^\# x\parallel_{\It}\leq \parallel x\parallel
\ \parallel T\parallel_{\It} ,
$$
is easy. This is because, for each $x\in X$, the operator
$T^\#x$ is the composition of $T$ with the bounded linear operator 
from
$Y$ to $X\hat \otimes_\epsilon Y$, which to each $y$ in $Y$ gives the
element $x\otimes y$. 

\medskip

If $i:\ Z\longrightarrow Z^{**}$ denotes the
isometric embedding of $Z$ into $Z^{**}$, it induces a bounded
linear operator $\hat i:\ I(Y,Z)\longrightarrow
I(Y,Z^{**})$ given by
$\hat
i(U)=i\circ U$ for all $U\in I(Y,Z)$.  It is immediate that $\hat i$ 
is
an isometry.  We will now show that the operator $\hat i\circ T^\#:\
X\longrightarrow I(Y,Z^{**})$\ is integral. It is well known (see 
[3,~p.~237])
that the space $I(Y,Z^{**})$ is isometrically isomorphic to the dual 
space
$(Y\hat \otimes_\epsilon Z^*)^*$. Thus to show that $\hat i\circ 
T^\#:\
X\longrightarrow (Y\hat \otimes_\epsilon Z^*)^*$ is an integral 
operator, we
need to show that it induces an element of $\left(X\hat 
\otimes_\epsilon
(Y\hat \otimes_\epsilon Z^*)\right)^*$.  For this, it is enough to 
note
that, by our discussion concerning the isometry of $(X\hat
\otimes_\epsilon Y)\hat \otimes_\epsilon Z^*$ and $X\hat
\otimes_\epsilon (Y\hat \otimes_\epsilon Z^*)$, that
$$
\parallel\hat i\circ 
T^\#\parallel_{\It}=\sup\left\{|\sum\limits^n_{i=1}\hat i\circ
T^\#x_i, y_i\otimes z^*_i|:\parallel\sum\limits^n_{i=1}x_i\otimes 
y_i\otimes
z^*_i\parallel_\epsilon \leq 1\right\} .\tag **
$$
But for each $x\in X$, $y\in Y$ and $z^*\in Z^*$, we have
$$
\langle \hat i\circ T^\#x, y\otimes z^*\rangle =\langle T (x\otimes 
y),
z^*\rangle .
$$
Hence, from (*) and (**), it follows that
$$
\parallel\hat i\circ T\parallel_{\It}=\parallel T\parallel_{\It}.
$$

Thus we have shown that (i) $\Rightarrow$ (ii). The proof of (ii) 
$\Rightarrow$
(i) follows in a similar way.  If $\hat i\circ T^\#:\
X\longrightarrow I(Y,Z^{**})$ is an integral operator, then one can 
show that
$i\circ T:\ X\hat \otimes_\epsilon Y\longrightarrow Z^{**}$ is 
integral, which in
turn implies that $T$ itself is integral (see [3,~p.~233]).
\medskip
Finally, the last assertion follows easily, since if $T^\#:\
X\longrightarrow I(Y,Z)$ is integral, then $\hat i\circ T$ is integral
(see [3,~p.~232]).
\endproof

Since the mapping $\hat i:\ I(Y,Z)\longrightarrow I(Y,Z^{**})$ is an
isometry, Proposition 1 coupled with Theorem 2 implies that, if $T:\
X\hat\otimes_\epsilon Y\longrightarrow Z$ is an integral operator, 
then
$T^\#:\ X\longrightarrow I(Y,Z)$ is 1-summing.  This result can be 
shown
directly from the definitions.  In what follows we shall present a 
sketch of
that alternative approach. 

\medskip

\noindent {\bf Theorem 3}\ \ Let $X$, $Y$ and $Z$ be Banach spaces, 
and let
$T:\ X\hat \otimes_\epsilon Y\longrightarrow Z$ be a bounded linear
operator.  If $T$ is integral, then $T^\#:\ X\longrightarrow I(Y,Z)$ 
is
1-summing.  If in addition $X$ is a {\it \$}$_\infty$-space, then $T:\
X\hat \otimes_\epsilon Y\longrightarrow Z$ is integral if and only if
$T^\#:\ X\longrightarrow I(Y,Z)$ is integral.

\medskip
\noindent {\bf Proof:}\ \ First, we will show that, if $T:\ X\hat
\otimes_\epsilon Y\longrightarrow Z$ is an integral operator, then 
$T^\#$ is in
$\prod_1\left(X,I(Y,Z)\right)$ with $\pi_1(T^\#)\leq \parallel 
T\parallel_{\text
{int}}$.   Let $x_1, x_2,\ldots, x_n$ be in
$X$, and fix
$\epsilon>0$.  For each $i\leq n$, there exists $n_i\in \Bbb N$,
$\left(y_{ij}\right)_{j\leq n_i}$ in $Y$, and $(z^*_{ij})_{j\leq n_i}$ 
in
$Z^*$, such that $\parallel \sum\limits^{ni}_{j=1}y_{ij}\otimes
z^*_{ij}\parallel_\epsilon \leq 1$, and
$$
\parallel T^\#x_i\parallel_{\text {int}}\leq
\sum\limits^{n_i}_{j=1}z^*_{ij}\left(T(x_i\otimes
y_{ij})\right)+\dfrac \epsilon {2^i}.
$$
Since $T$ is an integral operator, and
$$
\parallel \sum\limits^n_{i=1}\sum\limits^{n_i}_{j=1}x_i\otimes
y_{ij}\otimes z^*_{ij}\parallel_\epsilon \leq \sup
\left\{\sum\limits^n_{i=1}|x^*(x_i)|:\ \parallel x^*\parallel \leq 1,
x^*\in X^*\right\},
$$
it follows that
$$
\sum\limits^n_{i=1}\sum\limits^{n_i}_{j=1} z^*_{ij}
\left(T(x_i\otimes y_{ij})\right)\leq
\parallel T\parallel_{\text {int
}}\sup\left\{\sum\limits^n_{i=1}|x^*(x_i)|:\ \parallel x^*\parallel 
\leq
1,\ x^*\in X^*\right\}.
$$
Therefore
$$
\sum\limits^n_{i=1}\parallel T^\#x_i\parallel_{\text {int }}\leq
\parallel T\parallel_{\text {int
}}\sup\left\{\sum\limits^n_{i=1}|x^*(x_i)|:\ x^*\in X^*, \parallel
x^*\parallel \leq 1\right\}+ \epsilon.
$$

\medskip

Now, if in addition
$X$ is a {\it \$}$_\infty$-space, then by [12,~p.~477], the
operator $T^\#$ is indeed integral.
\endproof

\medskip

\noindent {\bf Remark 4}\ \ If $X=C(\Omega)$ is a space of
continuous functions defined on a compact Hausdorff space $\Omega$, 
one
can deduce a similar result to Theorem 3 from the main result of [13].

\medskip

Our next result extends a result of [16] to {\it \$}$_\infty$-spaces,
where it was shown that whenever
$X=C(\Omega)$, a space of all continuous functions on a compact
Hausdorff space $\Omega$, then a bounded linear operator $T:\
C(\Omega)\hat \otimes_\epsilon Y\longrightarrow Z$ is 1-summing if and 
only
if $T^\#:\ C(\Omega)\longrightarrow \prod_1(Y,Z)$ is 1-summing.  This 
also
extends a result of [14] where similar conclusions were shown to be 
true
for $X=A(K)$, a space of continuous affine functions on a Choquet
simplex $K$ (see [2]).

\medskip

We note that one implication follows with no restriction on $X$.
If $X$, $Y$ and $Z$ are Banach spaces, and $T:\ X\hat \otimes_\epsilon
Y\longrightarrow Z$ is a 1-summing operator, then $T^\#$ takes its
values in $\prod_1(Y,Z)$.  This follows from the fact that for each
$x\in X$, the operator $T^\#x$ is the composition of $T$ with the
bounded linear operator from $Y$ into $X\hat \otimes_\epsilon Y$ which
to each $y$ in $Y$ gives the element $x\otimes y$ in $X\hat
\otimes_\epsilon Y$, and hence
$$
\pi_1(T^\#x)\leq \parallel x\parallel \pi_1(T).
$$
Moreover, one can proceed as in [16] to show that $T^\#:\
X\longrightarrow \prod_1(Y,Z)$ is 1-summing.

\medskip

\noindent {\bf Theorem 5}\ \ If $X$ is a {\it \$}$_\infty$ space, then
for any Banach spaces $Y$ and $Z$, a bounded linear operator $T:\ 
X\hat
\otimes_\epsilon Y\longrightarrow Z$ is 1-summing if and only if 
$T^\#:\
X\longrightarrow \prod_1(Y,Z)$ is 1-summing.
\medskip
\noindent {\bf Proof:}\ \ Let $T:\ X\hat \otimes_\epsilon 
Y\longrightarrow
Z$ be such that $T^\#:\ X\longrightarrow \prod_1(Y,Z)$ is 1-summing.  
Since
$X$ is a {\it \$}$_\infty$-space, it follows from [14,~p.~477] that
$T^\#:\ X\longrightarrow \prod_1(Y,Z)$ is an integral operator. Let
$\varphi$ denote the isometric embedding of $Z$ into
$C\left(B(Z^*)\right)$, the space of all continuous scaler functions 
on
the unit ball $B(Z^*)$ of $Z^*$ with its weak$^*$-topology.  This 
induces
an isometry
$$
\hat \varphi:\ {\tsize\prod\nolimits}_1(Y,Z)\longrightarrow
{\tsize\prod\nolimits}_1\left((Y,C(B(Z^*))\right),
$$
$$
\hat \varphi (U)=\varphi\circ U\qquad \text { for all }U\in
{\tsize \prod\nolimits}_1(Y,Z).
$$
Now, it follows from [15,~p.~301], that 
$\prod_1\left(Y,C(B(Z^*))\right)$
is isometric to $I\left(Y, C(B(Z^*))\right)$. Hence we may assume that
$\hat \varphi \circ T^\#:\ X\longrightarrow I\left(Y,C(B(Z^*))\right)$
is an integral operator.  Moreover, it is easy to check that
$(\varphi\circ T)^\#=\hat \varphi \circ T^\#$.  By Theorem 2 the
operator $\varphi \circ T:\ X\hat \otimes_\epsilon Y\longrightarrow
C(B(Z^*))$ is an integral operator, and hence $T$ is in
$\prod_1\left(X\hat \otimes_\epsilon Y, Z\right)$ by Proposition 1.  
\endproof

In the following
section we shall, among other things, exhibit an example that 
illustrates that it
is crucial for the space $X$ to be a {\it \$}$_\infty$-space if the 
conclusion of
Theorem 5 is to be valid. 

\bigskip

\noindent {\bf III\ \ 2-summing Operators and some Counter-examples.}

In this section we shall study the behavior of 2-summing operators on
injective tensor product spaces.  As we shall soon see, the behavior 
of
such operators when $p=2$ is quite different from when $p=1$.  For
instance, unlike the case $p=1$, the {\it \$}$_\infty$-spaces don't 
seem to play
any particular role.  In fact, we shall exhibit operators $T$ on 
$C[0,1]\hat
\otimes_\epsilon \ell_2$ which are not 2-summing, yet their
corresponding operators $T^\#$ are.  We will also give other
interesting examples that answer some other natural questions.

\medskip

We will present the next theorem for $p=2$, but the same result is 
true for any
$1\leq p<\infty$, with only minor changes.

\medskip

\noindent {\bf Theorem 6}\ \ Let $X,Y$ and $Z$ be Banach spaces.  If
$T:\ X\hat \otimes_\epsilon Y\longrightarrow Z$ is a 2-summing 
operator,
then $T^\#:\ X\longrightarrow \prod_2 (Y,Z)$ is a 2-summing operator.

\medskip

\noindent {\bf Proof:}\ \ If $T:\ X\hat \otimes_\epsilon
Y\longrightarrow Z$ is 2-summing, then using the same kind of 
arguments that we
have given above, it can easily be shown that for each $x\in X$, that 
$T^\#x\in
\prod_2(Y,Z)$, with $\pi_2(T^\#x)\leq \pi_2(T)\parallel x\parallel$.

Now we will show that $T^\#:\ X\longrightarrow \prod_2 (Y,Z)$\ is 
2-summing.
Let $(x_n)$ be in $X$ such that $\sum\limits_n|x^*(x_n)|^2<\infty$\ 
for
each $x^*$ in $X^*$.  Fix $\epsilon>0$.  For each $n\geq 1$, let
$(y_{nm})$ be a sequence in $Y$ such that
$$
\sup\left\{\Bigg(\sum\limits^\infty_{m=1}|y^*(y_{nm})|^2\Bigg)^{1/2}:\
\parallel y^*\parallel \leq 1, y^*\in Y^*\right\}\leq 1 ,
$$
and 
$$
\pi_2\left(T^\#x_n\right)\leq \left(\sum\limits^\infty_{m=1}\parallel
T(x_n\otimes y_{nm})\parallel^2\right)^{1/2}+\dfrac \epsilon {2^n}.
$$
Then
$$
\left[\pi_2\left(T^\#x_n\right)\right]^2\leq
\sum\limits^\infty_{m=1}\parallel T\left(x_n\otimes
y_{nm}\right)\parallel^2+\dfrac \epsilon
{2^{n-1}}\left(\sum\limits^\infty_{m=1}\parallel T\left(x_n\otimes
y_{nm}\right)\parallel^2\right)^{1/2}+\dfrac {\epsilon^2}{2^{2n}} 
.\qquad
$$
Now, consider the sequence $(x_n\otimes y_{nm})$ in $X\hat 
\otimes_\epsilon
Y$.  For each $\xi \in \left(X\hat \otimes_\epsilon Y\right)^*\simeq
I(X,Y^*)$ we have that
$$
\aligned
\sum\limits_{m,n}\left|\xi(x_n)(y_{nm})\right|^2&=\sum\limits^\infty_{n
=1}
\sum\limits^\infty_{m=1}|\xi (x_n)(y_{nm})|^2\\
&\leq \sum\limits^\infty_{n=1}\parallel \xi (x_n)\parallel ^2.
\endaligned
$$
Since $\xi \in I(X,Y^*)$, it follows that $\xi\in \prod_2(X,Y^*)$, and 
so
$$
\sum\limits^\infty_{n=1}\parallel \xi (x_n)\parallel^2<\infty.
$$
Hence we have shown that for all $\xi \in (X\hat \otimes_\epsilon 
Y)^*$,
$$
\sum\limits_{m,n} \left|\xi (x_n)(y_{nm})\right|^2<\infty.
$$
Since $T\in
\prod_2\left(X\hat \otimes_\epsilon Y, Z\right)$, we have that
$$
\sum\limits_{m,n}\parallel T(x_n\otimes y_{nm})\parallel^2<\infty,
$$
and therefore
$$
\sum\limits_n\left[\pi_2\left(T{^\#}\!\!x_n\right)\right]^2<\infty.
$$
\endproof

\noindent {\bf Remark 7}\ \ The above result extends a result of [1], 
where
it was shown that if
$T:\ X\hat
\otimes_\epsilon Y\longrightarrow Z$ is $p$-summing for $1\leq
p<\infty$, then $T^\#:\ X\longrightarrow$ {\it \$}$(Y,Z)$ is
$p$-summing.

\medskip

Now we shall give the example that we promised at the end of 
section~II.

\medskip

\noindent {\bf Theorem 8}\ \ There exists a bounded linear operator 
$T:\
\ell_2\hat\otimes_\epsilon \ell_2\longrightarrow \ell_2$ such that
$T$ is not 1-summing, yet $T^\#:\ \ell_2\longrightarrow \pi_1(\ell_2,
\ell_2)$ is 1-summing.

\medskip

\noindent{\bf Proof:}\ \ First, we note the well known fact that 
$\ell_2\hat
\otimes_\epsilon \ell_2 = \slK(\ell_2, \ell_2)$, the space of all 
compact
operators from
$\ell_2$ to $\ell_2$. Now we define $T$\ as the composition of two 
operators. 

\medskip

Let
$P:\ \slK(\ell_2, \ell_2)\longrightarrow c_0$ be the operator defined 
so that for
each $K\in \slK(\ell_2, \ell_2)$,
$$
P(K)=\left(K(e_n)(e_n)\right),
$$
where $(e_n)$ is the standard basis of $\ell_2$.  It is well known
[10, p.145] that the sequence $(e_n\otimes e_n)$ in $\ell_2\hat
\otimes_\epsilon \ell_2$ is equivalent to the $c_0$-basis, and that 
the
operator $P$ defines a bounded
linear projection of $\slK(\ell_2, \ell_2)$ onto $c_0$.  

\medskip

Let
$S:\ c_0\longrightarrow \ell_2$ be the bounded linear operator such 
that
for each $(\alpha_n)\in c_0$
$$
S(\alpha_n)=\left(\dfrac {\alpha_n}n\right).
$$
It is easily checked [7,~p.~39] that $S$ is a 2-summing operator that 
is
not 1-summing.  

\medskip

Now we define $T:\ \slK(\ell_2, \ell_2)\longrightarrow
\ell_2$ to be $T=S\circ P$. Thus $T$ is 2-summing but not
1-summing.  It follows from Theorem 6 that the induced operator
$T^\#:\ell_2\longrightarrow \prod_2(\ell_2, \ell_2)$ is 2-summing.  
Since
$\ell_2$ is of cotype 2, it follows from [10,~p.~62], that for any 
Banach
space $E$, we have $\prod_2(\ell_2, E)=\prod_1(\ell_2, E)$, and that 
there exists
a constant $C>0$ such that for all $U\in \prod_2(\ell_2, E)$\ we have
$$
\pi_1(U)\leq C\pi_2(U).
$$
This implies that $T^\#$ is 1-summing as an operator
taking its values in $\prod_1(\ell_2, \ell_2)$.
\endproof

\medskip

\noindent {\bf Remark 9}\ \ 
We do not need to use Theorem~6 to show that $T^\#$\ is 1-summing in 
the
example above. Instead, we can use the following argument. First note 
that $T^\#$
factors as follows:
$$
\matrix \ell_2&@>{T^\#}>>&\pi_1(\ell_2, \ell_2)\\
@VAVV\\
\ell_2&\nearrow B\endmatrix
$$
Here $A:\ \ell_2\rightarrow \ell_2$ is the 1-summing operator defined 
by 
$$
A(\alpha_n)=\left(\dfrac {\alpha_n}n\right),
$$
for each $(\alpha_n)\in \ell_2$,
and $B:\ \ell_2\longrightarrow \pi_1(\ell_2, \ell_2)$ is the natural
embedding of $\ell_2$ into the space $\pi_1(\ell_2, \ell_2)$\ defined 
by
$$
B(\beta_n)(\gamma_n)=(\beta_n\gamma_n)
$$
for each $(\beta_n)$, $(\gamma_n)\in \ell_2$.

\bigskip

Now we will give two examples concerning the case when $p>1$. We will 
show that
we do not have a converse to Theorem~8, even when the underlying space 
$X$
is a {\it \$}$_\infty$-space.  

\medskip

First, let us fix some notation.
In what follows we shall denote the space $\ell_p(\Z)$\ by $\ell_p$, 
and
call its standard basis $\{e_n:n\in \Z\}$.  Thus if
$x=\left(x(n)\right)\in \ell_p$, then $x(n)=\langle x, e_n\rangle$, 
and
$$
\parallel x\parallel_{\ell_p}=\left(\sum\limits^\infty_{n=1}|\langle 
x,
e_n\rangle|^p\rangle\right)^{\frac 1p}.
$$
For $f\in L_p [0,1]$, we let
$$
\parallel f\parallel_{L_p}=\left(\int^1_0|f(t)|^pdt\right)^{\frac 1p}.
$$
If $\Omega$ is a compact Hausdorff space, and $Y$ is a Banach space, 
then
$C(\Omega,Y)=C(\Omega)\hat \otimes_\epsilon Y$ will denote the Banach
space of continuous $Y$-valued functions on $\Omega$ under the 
supremum
norm.  

\medskip

We recall that since $\ell_2$ is of cotype 2, we have that
$\prod_2(\ell_2,\ell_2)=\prod_1(\ell_2, \ell_2)$. We also recall that, 
if
$u=\sum\limits^\infty_{n=1}\alpha_ne_n\otimes e_n$ is a diagonal
operator in $\prod_2(\ell_2, \ell_2)$, then
$$
\pi_2(u)=\left(\sum\limits^\infty_{n=1}|\alpha_n|^2\right)^{\frac
12}=\text { the Hilbert-Schmidt norm of $u$.}
$$

\bigskip

\noindent {\bf Theorem 10}\ \ For each $1<p<\infty$, there is a
bounded linear operator $T:\ C([0,1], \ell_2)\rightarrow
\ell_2$
that is not $p$-summing, but such that $T^\#:\ C[0,1]\longrightarrow
\Pi_p(\ell_2, \ell_2)$ is $p$-summing.

\bigskip

\noindent {\bf Proof:}\ \ We present the proof for $p\le2$.
The case where $p>2$\ follows by the same argument. For each $n\in 
\Z$, let
$\epsilon_n(t):\ [0,1]\rightarrow \C,\ \epsilon_n(t)=e^{2\pi \It}$ 
denote the
standard trigonometric basis of $L_2[0,1]$.  If $f\in L_1[0,1]$, let 
$\hat
f(n)=\int^1_0f(t)\epsilon_n(t)dt$ denote the usual Fourier coefficient 
of $f$. 
For each $\lambda=(\lambda_n)$, where $|\lambda_n|\leq 1$ for all 
$n\in \Z$,
define the operator $$
T_\lambda:\ C\left([0,1], \ell_2\right)\longrightarrow \ell_2
$$
such that for $\varphi\in C\left([0,1], \ell_2\right)$\ we have
$$
T_\lambda\varphi=\left(\lambda_n\ \langle \hat \varphi (n), 
e_n\rangle\
\right).
$$
Here $\hat \varphi (n)=\text { Bochner 
--}\int^1_0\varphi(t)\epsilon_n(t)dt$. 

\medskip

The
operator $T_\lambda$ is a bounded
linear operator, with $\parallel T_\lambda \varphi 
\parallel_{\ell_2}\leq \parallel \varphi
\parallel$.  To see this, note that for $\varphi\in C\left([0,1],
\ell_2\right)$\ we have
$$
\aligned
\parallel T_\lambda \varphi \parallel^2_{\ell_2}&=\sum\limits_n
|\lambda_n|^2|\,\langle \hat \varphi (n), e_n\rangle\,|^2\\
&\leq \sum\limits_n|\,\langle \hat \varphi (n), e_n\rangle\, |^2\\
&\leq \sum\limits_n\int^1_0|\,\langle \varphi (t), e_n\rangle\, 
|^2dt\\
&=\int^1_0\parallel \varphi (t)\parallel^2_{\ell_2}dt\\
&\leq \sup\limits_t\parallel \varphi 
(t)\parallel^2_{\ell_2}.\endaligned
$$
Now, note that if $f\in C\left([0,1]\right)$, and $x\in \ell_2$, then
$$
T_\lambda (f\otimes x)=\left(\lambda_n \hat f(n)\langle x,
e_n\rangle\right),
$$
and hence the operator $T_\lambda^\#:\ C[0,1]\rightarrow ${\it 
\$}$(\ell_2,
\ell_2)$ is such that
$$
T^{\#}_\lambda f (x)=\left(\lambda_n\hat f(n) \langle x, e_n \rangle 
\right).
$$
Thus
$$
\pi_2(T^\#_\lambda f)=\left(\sum\limits_n |\lambda_n|^2|\hat
f(n)|^2\right)^{\frac 12} .$$
Hence, by H\"older's inequality,
$$
\pi_2(T^\#_\lambda f)
\leq \parallel (\lambda_n)\parallel_{\ell_r}\parallel (\hat f
(n))\parallel_{\ell_q},
$$
where $\dfrac 1r+\dfrac 1q=\dfrac 12$.
By the Hausdorff-Young inequality, we have that
$$
\parallel (\hat f (n))\parallel_{\ell_q}\leq \parallel
f\parallel_{L_p},
$$
where $1\leq p\leq 2$ and $\dfrac 1p+\dfrac 1q=1$.  Thus
$$
\pi_2(T^\#_\lambda f)\leq \parallel (\lambda_n)\parallel_{\ell_r}\ 
\parallel
f\parallel_{L_p},
$$
for $1\leq p\leq 2,\ 2\leq r\leq \infty$ and $\dfrac 1p=\dfrac
1r+\dfrac 12$. 
This shows that if $\parallel
(\lambda_n)\parallel_{\ell_r}<\infty$, then
\item {(1)} $T^\#_\lambda\left(C[0,1]\right)\subseteq \pi_2 (\ell_2,
\ell_2)=\pi_p(\ell_2, \ell_2)$;
\item {(2)} $T^\#_\lambda:\ C[0,1]\longrightarrow
\pi_p(\ell_2, \ell_2)$ is $p$-summing.

\medskip

Now, let $U\subset C\left([0,1], \ell_2\right)$ be the
closed linear span of $\{\epsilon_i\otimes e_i,\ a_i\in
\Z\}$. Then $U$ is isometrically isomorphic to $\ell_2$. This is 
because
$$
\aligned
\parallel \sum\limits_i\mu_i\epsilon_i\otimes e_i\parallel
&=\sup\limits_{t\in [0,1]}\parallel
\left(\mu_n\epsilon_n(t)\right)\parallel_{\ell_2}\\
&=\parallel \left(\mu_i\epsilon_i
(t_0)\right)\parallel_{\ell_2},\endaligned
$$
for some $t_0 \in [0,1]$, and hence
$$\parallel \sum\limits_i\mu_i\epsilon_i\otimes
e_i\parallel=\left(\sum\limits_i|\mu_i|^2\right)^{\frac 12}.
$$
Moreover
$$
T_\lambda(\epsilon_i\otimes e_i)=\lambda_ie_i\qquad \text { for all
}i\in \Z,
$$
Therefore, we have the following commuting diagram
$$
\matrix
U&@>{T_{\lambda |U}}>>&\ell_2\\
@VQVV&\nearrow {S_\lambda}&\\
\ell_2\endmatrix
$$
where $Q:\ U\rightarrow \ell_2$ is the isomorphism from $U$ onto
$\ell_2$ such that $Q(\epsilon_n\otimes e_n)=e_n$ for all $n\in \Z$, 
and
$S_\lambda:\ \ell_2\longrightarrow \ell_2$ is the operator given by
$S_\lambda(e_n)=\lambda_n e_n$.  So to show that $T_\lambda$ is not
$p$-summing, it is sufficient to show that one can pick
$\lambda=(\lambda_n)$ such that $S_\lambda$ is not $p$-summing. To do 
this, we
consider two cases. If $p=2$, we take $\lambda_n=1$ for all $n\in \Z$. 
 Then the
map $S_\lambda$ induced on $\ell_2$ is the identity map which is not
$s$-summing for any $s<\infty$.  If $1<p<2$, let $\lambda_n=\dfrac
1{|n+1|^{\frac 1r}\log |n+1|}$, so that $\parallel
(\lambda_n)\parallel_{\ell_r}<\infty$.  Then the map $S_\lambda:\
\ell_2\longrightarrow \ell_2$ is not $s$-summing for any $s<r$.  To
show this, we may assume, without loss of generality, that $s\geq
2$.  Let $x_n=e_n$ for all $n\geq 1$, and note that
$$\sup\limits_{x^*\in
B(\ell_2)}\left(\sum\limits_n|x^*(x_n)|^{s}\right)^{\frac 1{s}}\leq
\parallel x^*\parallel_{\ell_2}\leq 1,
$$
whilst
$$
\left(\sum\limits_n \parallel
\lambda_nx_n\parallel^{s}\right)^{\frac 1{s}}=\infty.
$$
\endproof

While the operators $T_\lambda$ in the previous example failed to be
$p$-summing, they were all (2,1)-summing.  This suggests the following
question:\ \ suppose $T:\ C\left([0,1], Y\right)\longrightarrow Z$ is 
a
bounded linear operator such that $T^\#:\ C[0,1]\longrightarrow
\prod_2(Y,Z)$ is 2-summing.  What can we say about $T$?  Is $T$\
$(2,1)$-summing?  The following example shows that $T$ can be very 
bad. 

\bigskip

\noindent {\bf Theorem 11}\ \ There exists a Banach space $Z$, and a 
bounded
linear operator \hfill\break
$T:\ C\left([0,1], \ell_1\right)\rightarrow Z$ such that
$T^\#:\ C[0,1]\rightarrow \prod_2(\ell_1, Z)$\ is $2$-summing, with 
the
property that, for any $N\in \N$, there exists a subspace $U$ of 
$C\left([0,1],
\ell_1\right)$ with $\dim U=N$, such that $T$ restricted to $U$ 
behaves like
the identity operator on $\ell^N_\infty$.  In particular $T$ is not
(2,1)-summing. 
\bigskip 
\noindent {\bf Proof:}\ \ If $X$ and $Y$ are Banach
spaces, we denote by $X \hat \otimes_\pi Y$\ the projective tensor 
product,
that is, the completion of the algebraic tensor product of $X$ and $Y$ 
under the
norm 
$$\parallel u\parallel_\pi=\inf\{\sum\limits^n_{i=1}\parallel x_i
\parallel \parallel y_i\parallel,\ u=\sum\limits^n_{i=1}x_i\otimes
y_i\}.
$$
It is well known that $(X\hat \otimes_\pi Y)^*$ is isometrically
isomorphic to the space {\it \$}$(X,Y^*)$ of all bounded linear
operators from $X$ to $Y^*$.  

\medskip

Let $Z=C\left([0,1],
\ell_1\right)+L_2[0,1]\hat \otimes_\pi \ell_2$\ be the Banach space 
with
the norm
$$
\parallel x\parallel_Z=\inf\{\parallel x'\parallel_\epsilon +\parallel
x''\parallel_\pi:\ x=x'+x''\},
$$
where $\parallel\ \parallel_\epsilon$ denotes the sup norm in
$C\left([0,1], \ell_1\right)$, and $\parallel\ \parallel_\pi$ denotes 
the
norm of the projective tensor product $L_2[0,1]\hat \otimes_\pi
\ell_2$.  Let
$$
T:\ C\left([0,1], \ell_1\right)\longrightarrow Z
$$
be the identity operator.  

\medskip

We first see that for each $f\in C[0,1]$, the operator $T^\#f:\
\ell_1\rightarrow Z$ is 2-summing with
$$
\pi_2(T^\#f)\leq \pi_2(I)\parallel T^\#f\parallel_{\text {\it
\$}(\ell_2,
Z)},
$$
where $I:\ell_1\longrightarrow \ell_2$ is the natural mapping.
This is because, for each $f\in C[0,1]$, and each $x\in
\ell_1$, we have that
$$
\parallel T(f\otimes x)\parallel \leq 
\parallel f\otimes x\parallel_{L_2\hat \otimes_\pi \ell_2} \leq
\parallel f\parallel_{L_2} \ 
\parallel x\parallel_{\ell_2}.
$$
To see that $T^\#:\ C[0,1]\longrightarrow
\prod_2(\ell_1, X)$\ is 2-summing, note that $\parallel 
T^\#f\parallel_{\text
{\it \$}(\ell_2, Z)}\leq \parallel f\parallel_{L_2}$, and hence if
$f_1, \ldots, f_n \in C[0,1]$, then
$$
\aligned
\left(\sum\limits^n_{k=1}\left[\pi_2(T^\#f_k)\right]^2\right)^{\frac
12}&\leq \pi_2(I) \left(\sum\limits^n_{k=1}\parallel
f_k\parallel^2_{L_2}\right)^{\frac 12}\\
&\leq \pi_2 (I) \pi_2(J) \sup\limits_{t\in [0,1]}\left\|
\left(\sum\limits^n_{K=1}|f_k(t)|^2\right)^{\frac
12}\right\|.\endaligned
$$
Here $J:\ C[0,1]\longrightarrow L_2[0,1]$ denotes the natural mapping.

\medskip

Now we define the space $U$, a
closed linear subspace of $C\left([0,1], \ell_1\right)$. Let 
$\{f_{ij}:\ 1\leq i, j\leq N\}$ be disjoint
functions in $C[0,1]$, for which $0\leq f_{ij} \leq 1$, $ \parallel
f_{ij}\parallel =1$, each $f_{ij}$ is supported in an interval of 
length
$\dfrac 1{N^2}$, and
$$
\int^1_0f_{ij}dt=\dfrac 1{2N^2}\text { and }\int^1_0f_{ij}^2dt=\dfrac
1{3N^2}.
$$
Let $\{e_{ij}:\ 1\leq i, j\leq N\}$ be distinct unit vectors in 
$\ell_1$. We let
$U=\{\sum\limits_{i,j}\lambda_if_{ij}\otimes e_{ij},\ \lambda_i\in 
\R\}$.

Now we consider $T$\ restricted to $U$.  If
$\sum\limits_{i,j}\lambda_i f_{ij}\otimes e_{ij}\in U$, then $$
\parallel\sum\limits_{i,j}\lambda_if_{ij}\otimes 
e_{ij}\parallel_\epsilon \leq
\sup\limits_i|\lambda_i|,
$$
and hence
$$\parallel \sum\limits_{i,j}\lambda_if_{ij}\otimes
e_{ij}\parallel_Z\leq \sup\limits_{i}|\lambda_i|.
$$
Let $y^*_i=N\sum\limits_jf_{ij}\otimes e_{ij}$, and set
$x=\sum\limits_{i,j}\lambda_if_{ij}\otimes e_{ij}$.  Then whenever
$x=x'+x''$, with $x'\in C\left([0,1], \ell_1\right)$ and $x''\in
L_2[0,1]\hat \otimes_\pi \ell_2$, we know that
$$
|y^*_i(x)|\leq |y^*_i(x')|+|y^*_i(x'')|.
$$
Hence
$$
\aligned
|y^*_i(x)|&\leq \parallel y^*_i\parallel_{C\left([0,1], 
\ell_1\right)^*}\
\parallel x'\parallel_\epsilon+\parallel y^*_i\parallel_{(L_2[0,1]\hat
\otimes_\pi \ell_2)^*}\
\parallel x''\parallel_\pi.\endaligned
$$
But
$$
\aligned
\parallel y^*_i\parallel_{C\left([0,1],
\ell_1\right)^*}&=N\sum\limits^N_{i=1}\int_{\supp
f_{ij}}|f_{ij}|dt\\
&=N\cdot \dfrac N{2N^2}=\dfrac 12,\endaligned
$$
and, since $\left(L_2[0,1]\hat \otimes_\pi \ell_2\right)^*$
is isometric to {\it \$}$\left(L_2[0,1], \ell_2\right)$,
$$
\aligned
\parallel y^*_i\parallel_{(L_2[0,1]\hat \otimes_\pi \ell_2)^*}
&=\sup\left\{[\sum\limits^N_{j=1}(N\int^1_0f_{ij}gdt)^2]^{\frac 12}:\
\parallel g\parallel_{L_2}\leq 1\right\}\\
&\leq \sup \left\{N[\sum\limits^N_{j=1}\int^1_0f_{ij}^2dt\cdot
\int_{\supp f_{ij}}|g|^2dt]^{\frac 12}:\parallel g\parallel_{L_2}\leq
1\right\}\\
&=\dfrac 1{\sqrt 3}\left\{(\sum\limits^N_{j=1}\int_{\supp
f_{ij}}|g|^2dt)^{\frac 12}:\parallel g\parallel_2\leq 1\right\}\\
&=\dfrac 1{\sqrt 3}.\endaligned
$$
Therefore
$$
|y_i^*(x)| \leq \dfrac 12 \parallel x'\parallel_\epsilon +\dfrac 
1{\sqrt
3}\parallel x''\parallel_\pi,\leq \dfrac 1{\sqrt 3}\parallel 
x\parallel.
$$
However,
$$
\aligned
y^*_i(x)&=N\sum\limits^N_{j=1}\lambda_i\int^1_0f_{ij}^2dt\\
&=N^2\lambda_i\dfrac 1{3N^2}=\dfrac {\lambda_i}3.\endaligned
$$
Therefore
$$
\aligned
\parallel \sum\limits_{i,j}\lambda_if_{ij}\otimes 
e_{ij}\parallel_Z&\geq
\sqrt 3 \sup\limits_i |y^*_i(x)|\\
&\geq \dfrac 1{\sqrt 3}\sup |\lambda_i|.\endaligned
$$
Thus the space $U$ is isomorphic to $\ell^N_\infty$, and we have the
commuting diagram
$$
\CD
U @>{T_{|U}}>> T(U)\\
@VAVV          @AA{A^{-1}}A\\
\ell^N_\infty  @>{id_{\ell^N_\infty}}>> \ell^N_\infty\endCD
$$
where $A:\ U\rightarrow \ell^N_\infty$ is the isomorphism between $U$ 
and
$\ell^N_\infty$.
\endproof

\bigskip

\noindent{\bf IV Operators that factor through a Hilbert space}

It is well known that {\it \$}$(X,\ell_2)=\prod_2(X,
\ell_2)$\ whenever $X$\ is $C(K)$\ or $\ell_1$. One might ask whether 
this is
true when $X = C(K,\ell_1)$. Indeed one could ask the weaker question: 
if $T:\
C(K, \ell_1)\longrightarrow \ell_2$\ is bounded, does it follow that 
the
induced operator $T^\#$\ is 2-summing? We answer this question in the 
negative.

\bigskip

\noindent {\bf Theorem 12}\ \ There is a compact Hausdorff space $K$ 
and
a bounded linear operator $T:\ C(K, \ell_1)\longrightarrow \ell_2$ for
which $T^\#:\ C(K)\longrightarrow \prod_1(\ell_1, \ell_2)$ is not
2-summing.

\medskip

\noindent {\bf Proof:}\ \  First, we show that there is a compact 
Hausdorff
space $K$, and an operator $R:\ C(K)\longrightarrow \ell_\infty$ that 
is
(2,1)-summing but not 2-summing. To see this, let $K=[0,1]$, and 
consider the
natural embedding $C[0,1] \longrightarrow L_{2,1}[0,1]$, where 
$L_{2,1}[0,1]$\
is the Lorentz space on $[0,1]$\ with the Lebesque measure (see [6]). 
By
[11], it follows that this map is (2,1)-summing. To show that this map
is not 2-summing, we argue in a similar fashion to [8]. For $n\in\N$,
consider the functions $e_i(t) = f(t+\frac 1i \bmod 1)$\ ($1\le i\le 
n$), where
$ f(t) = \frac 1{\sqrt t}$\ if $t\ge \frac 1n$\ and $\sqrt n$\
otherwise. Then it is an easy matter to verify that for some constant
$C>0$, $$ \left( \sum_{i=1}^n |e^*(e_i)|^2 \right)^{\frac12} \le C 
\sqrt{\log n}
$$ for every $e^*$\ in the unit ball of $C[0,1]^*$, whereas
$$ \left( \sum_{i=1}^n \| e_i \|_{L_{2,1}[0,1]}^2 \right)^{\frac12}
\ge C^{-1} \log n .
$$ 
Finally, since $L_{2,1}[0,1]$\ is separable, it embeds isometrically
into $\ell_\infty$.

\medskip

Define $T:\ C(K, \ell_1)\rightarrow
\ell_2$ as follows:\ \ for $\varphi=(f_n)\in C(K,\ell_1)$, let
$$
T(f_n)=\sum\limits_n Rf_n (n)e_n.
$$
Then $T$ is bounded, for
$$
\aligned
\parallel
T(f_n)\parallel_2&=\left(\sum\limits_n|Rf_n(n)|^2\right)^{\frac 12}\\
&\leq \left(\sum\limits_n \parallel
Rf_n\parallel^2_{\ell_\infty}\right)^{\frac 12}\\
&\leq \pi_{2,1}(R)\sup\limits_{t\in 
K}\sum\limits_n|f_n(t)|.\endaligned
$$
Thus
$$
\parallel T\parallel \leq \pi_{2,1}(R).
$$
But $T^\#:\ C(K)\longrightarrow$ {\it \$}$(\ell_1, \ell_2)$ is not
2-summing, because for each $f\in C(K)$, the operator $T^\#f:\
\ell_1\longrightarrow \ell_2$ is the diagonal operator
$\sum\limits_nRf(n)e_n\otimes e_n$.  Hence the strong operator norm of
$T^\#f$ is
$$
\parallel T^\#f\parallel =\sup\limits_n|Rf(n)|=\parallel
Rf\parallel_{\ell_\infty}.
$$
Thus $T^\#:\ C(K)\longrightarrow$ {\it \$}$(\ell_1, \ell_2)$ is not
2-summing, because $R:\ C(K)\longrightarrow \ell_\infty$ is not 
2-summing.
\endproof

\bigskip
\noindent {\bf Discussions and concluding remarks}

\noindent{\bf Remark 13}\ \ 
Theorem 12 shows that if $X$ and $Y$ are Banach spaces such that {\it 
\$}$(X,
\ell_2)=\prod_2(X, \ell_2)$ and {\it \$}$(Y, \ell_2)=\prod_2(X,
\ell_2)$, then $X\hat \otimes_\epsilon Y$ need not share this 
property.
This observation could also be deduced from arguments presented in [4] 
(use
Example~3.5 and the proof of Proposition~3.6 to show that there is a 
bounded
operator $T:(\ell_1\oplus\ell_1\oplus\ldots\oplus\ell_1)_{\ell_\infty}
\longrightarrow \ell_2$\ that is not $p$-summing for any $p<\infty$).

\medskip

\noindent {\bf Remark 14}\ \ In the proof of Theorem 2 we showed that 
the
injective tensor product is an associative operation, that is, if $X, 
Y$ and $Z$
are Banach spaces, then $(X\hat\otimes_\epsilon Y)\hat 
\otimes_\epsilon Z$ is
isometrically isomorphic to $X\hat \otimes_\epsilon (Y\hat 
\otimes_\epsilon
Z)$.  It is not hard to see that the same is true for the projective 
tensor
product. However, we can conclude from Theorem~12 that what is known 
as the
$\gamma_2^*$-tensor product is not an associative operation.

\medskip

If $E$ and $F$ are
Banach spaces, and $T:\ E\longrightarrow F$ is a bounded linear 
operator,
following [10], we say that $T$ {\bf factors through a Hilbert space} 
if
there is a Hilbert space $H$, and operators $B:\ E\longrightarrow H$ 
and $A:\
H\longrightarrow F$ such that $T=A\circ B$.  We let 
$\gamma_2(T)=\inf\{\parallel
A\parallel\ \parallel B\parallel \}$, where the infimum runs over all 
possible
factorization of $T$, and denote the space of all operators $T:\
E\longrightarrow F$ that factor through a Hilbert space by 
$\Gamma_2(E,F)$.  It
is not hard to check that $\gamma_2$ defines a norm on 
$\Gamma_2(E,F)$, making
$\Gamma_2(E,F)$\ a Banach space. We define the $\gamma_2^*$-norm
$\parallel\quad \parallel_*$ on $E\otimes F$ (see [9] or [10])
in which the dual of $E\otimes F$ is identified with
$\Gamma_2(E,F^*)$, and
let $E\hat \otimes_{\gamma^*_2}F$ denote the completion of $(E\otimes 
F,
\parallel\quad \parallel_*)$.

\medskip

The operator $T:\ C(K)\hat
\otimes_{\gamma^*_2}\ell_1\longrightarrow \ell_2$\ exhibited in 
Theorem 12,
induces a bounded linear functional on $\left[(C(K)\hat
\otimes_{\gamma^*_2}\ell_1)\hat \otimes_{\gamma^*_2}\ell_2\right]^*$.  
Now we
see that if $C(K)\hat \otimes_{\gamma^*_2}(\ell_1\hat
\otimes_{\gamma^*_2}\ell_2)$ were isometrically isomorphic to
$(C(K)\hat \otimes_{\gamma^*_2}\ell_1)\hat
\otimes_{\gamma^*_2}\ell_2$, then the operator $T^\#:\
C(K)\rightarrow ${\it \$}$(\ell_1, \ell_2)$ would induce a bounded 
linear
functional on $\left[C(K)\hat \otimes_{\gamma^*_2}(\ell_1\hat
\otimes_{\gamma^*_2}\ell_2)\right]^*$, showing that $T^\#\in
\Gamma_2\left(C(K),\text {\it \$}(\ell_1, \ell_2)\right)$, implying
that $T^\#$ would be 2-summing [10,~p.~62].  This contradiction shows 
that $C(K)
\hat \otimes_{\gamma^*_2}(\ell_1\hat \otimes_{\gamma^*_2}\ell_2)$ and
$\left(C(K)\hat \otimes_{\gamma^*_2}\ell_1\right)\hat
\otimes_{\gamma^*_2}\ell_2$ cannot be isometrically isomorphic.

\medskip

Another example showing that the $\gamma_2^*$-tensor product is not
associative was given by Pisier (private communication).

\vfill\eject

\centerline {\bf Bibliography}

\item {[1]} R. Bilyeu, and P. Lewis, {\it Some Mapping Properties of
Representing Measures}, Ann. Math Pure Appl. CIX (1976) p. 273--287.

\item {[2]} G. Choquet, {\bf Lectures on Analysis}, Vol. {\bf II},
Benjamin, New York, (1969).

\item {[3]} J. Diestel, and J.J. Uhl Jr., {\bf Vector Measures}, Math
Surveys, {\bf 15}, AMS, Providence, RI (1977).

\item{[4]} T. Figiel, J. Lindenstrauss, and V. Milman, {\it The 
dimension of
almost spherical sections of convex bodies}, Acta Mathematica, {\bf 
139},
(1977), p.~53--94.

\item {[5]} A. Grothendick, {\bf Produits tensoriels topologiques et
espaces nucl\'eaires}, Mem. A.M.S. {\bf 16}, (1955).

\item {[6]} R.A. Hunt, {\it On $L(p,q)$\ spaces}, L'Enseignement 
Math.\ (2),
{\bf 12}, (1966), p.~249--275.

\item {[7]} G.J.O. Jameson, {\bf Summing and Nuclear Norms in Banach
Space Theory}, LMSST {\bf 8}, Cambridge University Press (1987).

\item {[8]} G.J.O. Jameson, {\it Relations between summing norms of 
mappings on
$\ell_\infty$}, Math.~Z, {\bf 194}, (1987), p.~89--94.

\item {[9]} S. Kwapien, {\it On operators factorizable through 
$L_p$-spaces},
Bull.\ Soc.\ Math.\ France, M\'em 31--32, (1972), p.~215--225.

\item {[10]} G. Pisier, {\bf Factorization of Linear Operators and
Geometry of Banach Spaces}, AMS CBMS {\bf 60}, Providence RI (1986).

\item {[11]} G. Pisier, {\it Factorization of operators through 
$L_{p\infty}$\
or $L_{p1}$\ and non commutative generalizations}, Math.\ Ann., {\bf 
276},
(1986), p.~105--136.

\item {[12]} J.R. Retherford, and C. Stegall, {\it Fully Nuclear and
Completely Nuclear Operators with applications to {\it \$}$_1$ and 
{\it
\$}$_\infty$ spaces}, T.A.M.S., {\bf 163}, (1972) p. 457--492.

\item {[13]} P. Saab, {\it Integral Operators on Spaces of Continuous
Vector Valued Functions}, Proc.\ Amer.\ Math.\ Soc.\ (to appear).

\item {[14]} B. Smith, {\it Some Bounded Linear Operators On the 
Spaces
$C(\Omega,E)$ and $A(K,E)$}, Ph.D. Dissertation, The University of
Missouri-Columbia, 1989.

\item {[15]} C. Stegall, {\it Characterization of Banach spaces whose
duals are $L_1$ spaces}, Is. J. of Math, {\bf 11}, (1972) p. 299--308.

\item {[16]} C. Swartz, {\it Absolutely summing and dominated 
operators
on spaces of vector-valued continuous functions}, T.A.M.S., {\bf 179},
(1973) p. 123--132.

\vfill
\centerline{\vbox{\halign{#\hfill\cr
University of Missouri\cr
Dept.\ of Math.\cr
Columbia, MO 65211\cr}}}

\eject

\bye